\documentclass[12pt]{amsart}
\usepackage{amstext,amsfonts,amssymb,amscd,amsbsy,amsmath,verbatim, mathrsfs, fullpage}
\usepackage[alphabetic,abbrev,lite,backrefs]{amsrefs} 
\usepackage{ifthen,tikz}
\usepackage{color}
\usepackage{amsthm}
\usepackage{latexsym}
\usepackage[all]{xy}
\usepackage{enumerate}
\usepackage{mathtools}

\numberwithin{equation}{section}

\newtheorem{lemma}{Lemma}[section]

\newtheorem{prop}[lemma]{Proposition}

\newtheorem{claim*}{Claim}
\newtheorem{thm}[lemma]{Theorem}

\theoremstyle{definition}
\newtheorem{defn}[lemma]{Definition}
\newtheorem{example}[lemma]{Example}

\newtheorem{algorithm}[lemma]{Algorithm}

\theoremstyle{remark}

\newcommand{\cC}{\mathcal{C}}

\newcommand{\m}{\mathfrak m}
\newcommand{\PP}{\mathbb P}

\newcommand{\im}{\operatorname{im}}

\newcommand{\Hom}{\operatorname{Hom}} 

\newcommand{\defi}[1]{\textsf{#1}} 

\newcommand{\beq}{\begin{displaymath}}
\newcommand{\eeq}{\end{displaymath}}

\def\nc{\newcommand}
\def\on{\operatorname}
\nc{\Q}{\mathbb{Q}}
\nc{\RR}{\mathbf{R}}
\nc{\LL}{\mathbf{L}}
\nc{\xra}{\xrightarrow}
\nc{\xla}{\xleftarrow}

\def\om{\omega}

\def\DM{\operatorname{DM}}

\def\th{\on{th}}

\nc{\into}{\hookrightarrow}
\nc{\onto}{\twoheadrightarrow}
\nc{\OO}{\mathcal{O}}
\nc{\Z}{\mathbb{Z}}
\nc{\cA}{\mathcal{A}}
\nc{\w}{\widehat}
\nc{\End}{\on{End}}
\nc{\res}{\frac{1}{x_0x_1}}
\nc{\tF}{\widetilde{F}}
\nc{\tG}{\widetilde{G}}
\nc{\tf}{\widetilde{f}}
\nc{\Com}{\on{Com}}

\nc{\G}{\mathbb{G}}
\nc{\cG}{\mathcal{G}}
\nc{\cE}{\mathcal E}
\nc{\cF}{\mathcal F}
\nc{\cR}{\mathcal R}
\nc{\cD}{\mathcal D}
\nc{\cB}{\mathcal B}
\nc{\cT}{\mathcal T}
\nc{\cL}{\mathcal L}

\nc{\bM}{\mathbf M}
\nc{\bN}{\mathbf N}
\nc{\U}{\mathbf U}
\nc{\BM}{\mathbf B \mathbf M}
\nc{\Dsg}{\on{D}_{\on{sg}}}
\nc{\fC}{\mathcal{C}}
\nc{\fG}{\mathcal{G}}
\nc{\N}{\mathbb{N}}

\nc{\del}{\partial}
\nc{\cone}{\on{cone}}
\nc{\D}{\on{D}_{\on{diff}}}
\nc{\DMb}{\on{D}^b_{\DM}}
\nc{\Db}{\on{D}^{\on{b}}}
\nc{\Kb}{\on{K}^{\on{b}}}
\nc{\fm}{\mathfrak{m}}
\nc{\Flag}{\on{Flag}}
\nc{\DMmin}{\DM_{\on{min}}}
\nc{\Ddiff}{\on{D}_{\on{diff}}}
\nc{\Dbdiff}{\on{D}^\on{b}_{\on{diff}}}
\nc{\wO}{\widehat{\OO}}
\nc{\wT}{\widehat{T}}
\nc{\from}{\leftarrow}
\nc{\wLL}{\widetilde{\LL}}
\nc{\augCech}{\widetilde{\cC}}
\nc{\Fold}{\on{Fold}}
\nc{\Ext}{\on{Ext}}
\nc{\FF}{\mathbf{F}}
\nc{\Comper}{\Com_{\on{per}}}
\nc{\Unfold}{\on{Unfold}}
\nc{\intHom}{\underline{\Hom}}
\nc{\Ex}{\on{Ex}}
\nc{\tg}{\widetilde{g}}

\nc{\B}{\mathcal{B}}
\nc{\K}{\mathcal{K}}
\nc{\kos}{\on{Kos}}
\nc{\Perf}{\on{Perf}}
\nc{\tR}{\widetilde{\cR}}
\nc{\X}{\mathcal{X}}
\nc{\Cl}{\on{Cl}}
\nc{\fU}{\mathcal{U}}
\nc{\bU}{\mathbf U}
\def\c{\colon}
\nc{\st}{\on{st}}

\nc{\coh}{\on{coh}}

\def\D{\mathcal{D}}

\nc{\tU}{\U}
\nc{\bC}{\mathbf{C}}
\nc{\aux}{\on{aux}}
\def\MR#1{}
\nc{\co}{\colon}
\def\c{\colon}
\def\e{\varepsilon}

\def\ce{\coloneqq}

\title{The multigraded BGG correspondence in Macaulay2}

\thanks{We gratefully acknowledge support from NSF grant DMS-2302373 and the Natural Sciences and
Engineering Research Council of Canada (NSERC)}

\author{Maya Banks}
\author{Michael K. Brown}
\author{Tara Gomes}
\author{Prashanth Sridhar}
\author{Eduardo Torres Davila}
\author{Sasha Zotine}

\newcommand{\Addresses}{{
	\vskip\baselineskip
  	\footnotesize
  	\noindent \textsc{Department of Mathematics and Statistics, Auburn University} \par\nopagebreak
	\noindent \textit{E-mail addresses:} \texttt{mkb0096@auburn.edu, prashanth.sridhar0218@gmail.com}
  	\vskip\baselineskip
  	\noindent \textsc{School of Mathematics, University of Minnesota} \par\nopagebreak
	\noindent \textit{E-mail addresses:} \texttt{gomes072@umn.edu, torre680@umn.edu}
 \vskip\baselineskip
  	\noindent \textsc{Department of Mathematics and Statistics, Queen's University} \par\nopagebreak
	\noindent \textit{E-mail address:} \texttt{zotinea@mcmaster.ca}
	\vskip\baselineskip
  	\noindent \textsc{Department of Mathematics, University of Wisconsin} \par\nopagebreak
	\noindent \textit{E-mail address:} \texttt{mayadb@uic.edu}
 
}}

\begin{document}
\maketitle

\begin{abstract}
We give an overview of a \verb|Macaulay2| package for computing with the multigraded BGG correspondence. This software builds on the package \verb|BGG| due to Abo-Decker-Eisenbud-Schreyer-Smith-Stillman, which concerns the standard graded BGG correspondence. In addition to implementing the multigraded BGG functors, this package includes an implementation of differential modules and their minimal free resolutions, and it contains a method for computing strongly linear strands of multigraded free resolutions.
\end{abstract}


\section{Introduction}

The Bernstein-Gel'fand-Gel'fand (BGG) correspondence is an equivalence of derived categories between a polynomial and exterior algebra \cite{BGG}. 
The BGG correspondence has been widely applied in commutative algebra and algebraic geometry; we refer the reader to \cite[Chapter 7]{geosyz} for a detailed discussion. Many computational applications of the BGG correspondence have been implemented in the \verb|Macaulay2| package \verb|BGG| due to Abo-Decker-Eisenbud-Schreyer-Smith-Stillman. 

The BGG correspondence admits a generalization involving \emph{multigraded} polynomial rings, i.e. polynomial rings graded by $\Z^t$ for some $t \ge 1$. Just as the classical BGG correspondence plays an important role in projective geometry, the multigraded BGG correspondence plays an analogous role in toric geometry; see e.g. \cite{baranovsky, BETate, BEstrands, BEsyzygies, BEpositivity, BS, EES}. The goal of this article is to describe the \verb|Macaulay2| \cite{M2} package \verb|MultigradedBGG|, which builds on the \verb|BGG| package by implementing the multigraded BGG correspondence. 

Let us give an overview of the paper. In the multigraded BGG correspondence, one must consider not only \emph{complexes} of exterior modules, but more general objects called \emph{differential} exterior modules. The package \verb|MultigradedBGG| therefore includes an implementation of differential modules and their free resolutions: we discuss this aspect of the package in~\S\ref{DM}. In \S\ref{sec:BGG}, we recall the definitions of the multigraded BGG functors and compute some examples using \verb|MultigradedBGG|. In \S\ref{sec:strands}, we describe a method for computing strongly linear strands of multigraded free resolutions, in the sense of \cite{BEstrands}.



\section*{Acknowledgments} We are grateful to Daniel Erman and Keller VandeBogert for generously contributing many key insights to the code in  \verb|MultigradedBGG| and to Greg Smith  for many helpful comments during the writing of this package. We thank the Institute for Mathematics and its Applications at the University of Minnesota for hosting the “Macaulay2 Workshop \& Mini-School” (funded by NSF grant DMS-2302476) on June 3 - 9, 2023, during which
most of this package was developed. Finally, we thank the referees for many helpful suggestions.

\section{Computing with differential modules}
\label{DM}

We begin by recalling some background on differential modules. Let $R = \bigoplus_{i \ge 0} R_i$ be a (possibly noncommutative) graded ring such that $R_0$ is local.\footnote{While we assume $R$ is nonnegatively graded, versions of the results we state here hold for nonpositively graded rings as well. Indeed, all one needs is that $R$ is graded by an abelian group $A$ such that (a) $R_0$ local, and (b) $R$ is \defi{positively $A$-graded}, in the sense of \cite[Definition A.1]{BETate}.} 
All $R$-modules are right modules throughout. We denote by $\m$ the maximal ideal of $R$ given by the sum of the maximal ideal of $R_0$ and the ideal $\bigoplus_{i > 0} R_i$.

\begin{defn}
Given $a \in \Z$, a \defi{degree $a$ differential module} is a graded $R$-module $D$ equipped with a degree $a$ endomorphism $\del_D$ such that $\del_D^2 = 0$. The \defi{homology} of $(D, \del_D)$ is given by $\ker(\del_D \c D \to D(a)) / \im(\del_D \c D(-a) \to D)$ and denoted $H(D, \del_D)$. A \defi{morphism} $f \co (D, \del_D) \to (D', \del_{D'})$ of degree $a$ differential $R$-modules is a degree 0 morphism $f \co D \to D'$ such that $\del_{D'} f = f \del_D$.
\end{defn}

\begin{example}
\label{deg2}
Let $R$ denote the standard graded polynomial ring $k[x, y]$, where $k$ is a field. Set $D = R^2$ and $\del_D = \begin{pmatrix} xy & -x^2 \\ y^2 & -xy \end{pmatrix}$. The pair $(D, \del_D)$ is a degree 2 differential $R$-module.
\end{example}

There are natural notions of quasi-isomorphism, mapping cone, homotopy, and contractibility in the setting of differential modules generalizing the usual ones for chain complexes: see \cite[\S 2.1]{BE1} for the definitions. Given a morphism $f : (D, \del_D) \to (D', \del_{D'})$ of degree $a$ differential $R$-modules, we let $\cone(f)$ denote its mapping cone.

In the theory of free resolutions of differential modules developed in \cite{BE1}, a key role is played by differential modules called free flags, a notion due to Avramov-Buchweitz-Iyengar:

\begin{defn}[\cite{ABI} \S 2]
A degree $a$ differential $R$-module $(F, \del)$ is called a \defi{free
flag} if there exists a decomposition $F = \bigoplus_{i \ge 0} F_i$ such that each $F_i$ is free, 
and $\del(F_i)\subseteq \bigoplus_{j < i} F_j$. That is, $\del$ is block upper triangular with respect to the decomposition $F = F_0 \oplus F_1 \oplus \cdots$.
\end{defn}

\begin{example}
\label{ex:flag}
Let $R$ be as in Example~\ref{deg2}. An example of a degree 2 free flag is given by: 
$$
G = R(-1) \oplus R^2 \oplus R(1), \quad 
\del =
\begin{pmatrix}
0 & y & x & 1 \\ 
0 & 0 & 0 & x\\ 
0 & 0 & 0 & -y \\
0 & 0 & 0 & 0 \\
 \end{pmatrix}.
 $$
The differential module $(D, \del_D)$ from Example~\ref{deg2} is not isomorphic to a free flag~\cite[Example 5.6]{ABI}. But, there is a quasi-isomorphism $(G, \del_G) \xra{\simeq} (D, \del_D)$ \cite[Example 5.8]{BE1}.
\end{example}

We now recall the definition of a (minimal) free resolution of a differential module:

\begin{defn}
\label{def:minres}
Let $(D, \del_D)$ be a degree $a$ differential $R$-module. A \defi{free flag resolution} of $(D, \del_D)$ is a quasi-isomorphism $(G, \del_G) \xra{\simeq} (D, \del_D)$, where $(G, \del_G)$ is a free flag. Such a resolution is called a \defi{minimal free flag resolution} if $\del_G(G) \subseteq \m G$. A \defi{free resolution} of $(D, \del_D)$ is a quasi-isomorphism $\e \co (F, \del_F) \xra{\simeq} (D, \del_D)$ that factors as 
$$
\xymatrix{
(F, \del_F) \ar[rd]^-{\iota} \ar[rr]^-{\e} && (D, \del_D) \\
& (G, \del_G) \ar[ru]^-{\widetilde{\e}}, & \\
}
$$
where $\iota$ is a split injection of differential $R$-modules, and $\widetilde{\e}$ is a free flag resolution. We say $\e \co (F, \del_F) \xra{\simeq} (D, \del_D)$ is a \defi{minimal free resolution} if $\del_F(F) \subseteq \m F$.
\end{defn}

The following theorem gives sufficient conditions for the existence and uniqueness of minimal free resolutions of differential modules.

\begin{thm}[\cite{BE1} Theorem 1.2]
\label{EU}
Let $(D, \del_D)$ be a degree $a$ differential $R$-module. Assume its homology  $H(D, \del_D)$ is finitely generated. If either $a = 0$ or $H(D, \del_D)$ has finite projective dimension as an $R$-module, then $(D, \del_D)$ admits a minimal free resolution that is unique up to isomorphism.
\end{thm}


Any differential $R$-module admits a free flag resolution via the following algorithm, which is encoded as the function \verb|resDM| in the package \verb|MultigradedBGG|. This process closely resembles the usual algorithm for building a free resolution of an ordinary module.

\begin{algorithm}[\cite{BE1} Construction 2.8]
\label{alg:res}
Let $(D, \del_D)$ be a degree $a$ differential $R$-module. 
\begin{enumerate}
\item Choose a degree 0 map $\e_0 \co F_0 \to D$ such that $F_0$ is free, and $\e_0$ sends a basis of $F_0$ to cycles in $D$ that descend to a homogeneous minimal generating set for $H(D, \del_D)$. 
\item If $\cone(\e_0)$ is exact, we're done; $\e_0$ is a free flag resolution of $(D, \del_D)$. Otherwise, go back to Step (1) with $(D, \del_D)$ replaced with $\cone(\e_0)$.
\end{enumerate}
This algorithm may not terminate in finitely many steps, in which case passing to the colimit as $i \to \infty$ of the maps $\e_i \c F_i \to D$ yields a free flag resolution. 
\end{algorithm}

\begin{example}
\label{ex:res}
Let us use the \verb|resDM| function to build a free flag resolution of the differential module $D$ in Example~\ref{deg2}. 
\begin{footnotesize}
\begin{verbatim}

i1 : needsPackage "MultigradedBGG"
i2 : R = ZZ/101[x, y];
i3 : A = map(R^2, R^2, matrix{{x*y, -x^2}, {y^2, -x*y}}, Degree => 2);

\end{verbatim}
\end{footnotesize}

Applying the function \verb|differentialModule| to an object of type \verb|Matrix| yields an object of type \verb|DifferentialModule|, which is a \verb|Complex| with three properties: it is only nonzero in homological degrees $-1$, 0, and 1; its two differentials are identical; and the differentials compose to zero. 

\begin{footnotesize}

\begin{verbatim}

i4 : D = differentialModule(A);
      2      2      2
o4 : R  <-- R  <-- R
                     
    -1      0      1

o4 : DifferentialModule

\end{verbatim}
\end{footnotesize}

Let us now resolve $D$:
\begin{footnotesize}
\begin{verbatim}

i5 : F = resDM(D);
i6 : F.dd_0

o6 = {1}  | 0 y x 1  |
     {0}  | 0 0 0 x  |
     {0}  | 0 0 0 -y |
     {-1} | 0 0 0 0  |
               
\end{verbatim}
\end{footnotesize}
Observe that $F$ coincides with the free flag in Example~\ref{ex:flag}, confirming that this free flag is quasi-isomorphic to $D$.
\end{example}

As we see in Example~\ref{ex:res}, applying Algorithm~\ref{alg:res} does not always yield a \emph{minimal} free flag resolution. Indeed, such resolutions need not exist; for instance, the differential module $(D, \del_D)$ in Example~\ref{deg2} does not admit a minimal free flag resolution. To see this, observe that $(D, \del_D)$ is its own minimal free resolution \cite[Example 5.8]{BE1}; since it is not isomorphic to a free flag, the uniqueness of minimal free resolutions implies that it cannot admit a minimal free flag resolution. However, there is an algorithm for \emph{minimizing} any finitely generated free differential $R$-module $(F, \del_F)$, i.e. decomposing $(F, \del_F)$ as a direct sum $(M, \del_M) \oplus (C, \del_C)$, where $(M, \del_M)$ is minimal, and $(C, \del_C)$ is contractible: see \cite[Proposition 4.1]{BE1}. For instance, if $(F, \del_F)$ is a free flag resolution of a differential module $(D, \del_D)$, then $(M, \del_M)$ is the minimal free resolution of $(D, \del_D)$. This minimization algorithm is implemented as the function \verb|minimizeDM| in \verb|MultigradedBGG|. 

\begin{example}
\label{ex:minimize}
Let us minimize the free flag resolution obtained in Example~\ref{ex:res}. 
\begin{footnotesize}
\begin{verbatim}

i1 : needsPackage "MultigradedBGG"
i2 : R = ZZ/101[x, y];
i3 : A = map(R^2, R^2, matrix{{x*y, -x^2}, {y^2, -x*y}}, Degree => 2);
i4 : D = differentialModule(A);
i5 : F = resDM(D);
i6 : G = minimizeDM F;

\end{verbatim}
\end{footnotesize}
Checking the differential of the minimization $G$ of $F$ confirms that we recover the differential module $D$ that we started with (up to isomorphism); that is, $D$ is its own minimal free resolution, as observed in \cite[Example 5.8]{BE1}:
\begin{footnotesize}
\begin{verbatim}

i13 : G.dd_0

o13 = | -xy -x2 |
      | y2  xy  |

\end{verbatim}
\end{footnotesize}
\end{example}

While minimal free flag resolutions do not exist in general, we do have the following:

\begin{prop}
\label{prop:res}
Assume $R_0$ is a field. Let $(D, \del_D)$ be a degree 0 differential $R$-module with finitely generated homology. Choose $n$ such that $H(D)_n \ne 0$ and $H(D)_i = 0$ for $i < n$. There exists a minimal free flag resolution $(F, \del_F)$ of $(D, \del_D)$, and $(F, \del_F)$ may be equipped with a flag structure $F = \bigoplus_{i \ge 0} F_i$ such that each $F_i$ is finitely generated in degree $n + i$.
\end{prop}


A resolution as in Proposition~\ref{prop:res} is necessarily unique up to isomorphism, by Theorem~\ref{EU}. Our process for constructing such resolutions is nearly identical to Algorithm~\ref{alg:res}:

\begin{algorithm}
\label{alg:res2}
Since we may replace $D$ with $D(n)$, we may assume $n = 0$. We proceed via the following algorithm:
\begin{enumerate}
\item Let $y_1, \dots, y_t \in D$ be cycles that descend to a $R_0$-basis of $H(D)_0$. Set $F_0 = R^t$, and define a map $\e_0 \c F_0 \to D$ that sends the $i^{\th}$ basis element of $F_0$ to $y_i$.   
\item The map $\e_0$ determines an isomorphism on homology in degree 0. If $\cone(\e_0)$ is exact, we're done; $\e_0$ is a minimal free flag resolution of $(D, \del_D)$ with the desired properties. Otherwise, go back to Step (1) with $(D, \del_D)$ replaced with $\cone(\e_0)$, observing that $H(\cone(\e_0))_i = 0$ for $i < 1$. 
\end{enumerate}
As with Algorithm~\ref{alg:res}, this process need not terminate in finitely many steps, in which case passing to the colimit as $i \to \infty$ of the maps $\e_i \c F_i \to D$ gives a free flag resolution with the desired properties. The resolution $F$ is minimal since each $F_i$ is a direct sum of copies of $R(-i)$, and any degree 0 map $R(-i) \to \bigoplus_{0 \le j < i} R(-j)$ is minimal.
\end{algorithm}


Algorithm~\ref{alg:res2} is implemented in \verb|MultigradedBGG| as \verb|resMinFlag|. When a differential module $D$ admits a minimal free flag resolution, one may also obtain it by first resolving $D$ via Algorithm~\ref{alg:res} and then minimizing, but this is less efficient than Algorithm~\ref{alg:res2}.

\begin{example}
Let $R$ be the standard graded polynomial ring $k[x,y]$, where $k$ is a field, and let $D = k$, considered as a differential module with trivial differential. This differential module admits a minimal free flag resolution, namely the Koszul complex on $x$ and $y$ considered as a differential module; let us now use \verb|resMinFlag| to perform this simple calculation:
\begin{footnotesize}
\begin{verbatim}

i1 : needsPackage "MultigradedBGG"
i2 : R = ZZ/101[x, y];
i3 : k = coker vars R;
i3 : f = map(k, k, 0);
i4 : D = differentialModule(f);

\end{verbatim}
\end{footnotesize}
The function \verb|resMinFlag| takes two inputs: a differential module $D$ and an integer $t$ indicating how many times Algorithm~\ref{alg:res2} should be iterated. It follows from \cite[Theorem 3.2]{BE1} that, if
$$
d \ce \sup\{d - d' \text{ : } R(d) \text{ and } R(d') \text{ are summands of the minimal free resolution of $H(D)$}\},
$$
then when $d < \infty$, the differential module \verb|resMinFlag(D, d+1)| is the entire minimal free flag resolution of $D$. In our example, we have $d = 2$, and so we proceed as follows:
\begin{footnotesize}
\begin{verbatim}

i5 : F = resMinFlag(D, 3);
i6 : F.dd_0
           
o6 = {0} | 0 y x 0  |
     {1} | 0 0 0 x  |
     {1} | 0 0 0 -y |
     {2} | 0 0 0 0  |
                
\end{verbatim}
\end{footnotesize}
As expected, the output is the Koszul complex, considered as a differential module. 
\end{example}

\section{The multigraded BGG functors}
\label{sec:BGG}

Let $S = k[x_0, \dots, x_n]$, where $k$ is a field, and suppose $S$ is $A \ce \Z^t$-graded for some $t \ge 1$. For instance, $S$ could be the Cox ring of a smooth projective toric variety $X$ equipped with its $\Cl(X)$-grading, where $\Cl(X)$ is the divisor class group of $X$. Let $E$ denote the $A \oplus \Z$-graded exterior algebra  $\Lambda_k(e_0, \dots, e_n)$ with $\deg(e_i) = (-\deg(x_i); -1)$. We let $\Com(S)$ denote the category of complexes of $A$-graded $S$-modules and $\DM(E)$ the category of degree $(0 ; -1)$ differential right $E$-modules. Henceforth, all differential $E$-modules are assumed to have degree $(0; -1)$; i.e. they are assumed to be objects in $\DM(E)$.

As proven by \cite{HHW}, there is a multigraded analogue of the Bernstein-Gel'fand-Gel'fand (BGG) correspondence that gives an adjunction
$
\LL : \DM(E) \leftrightarrows \Com(S) : \RR.
$
We refer the reader to \cite[\S 2]{BETate} for a detailed introduction to the multigraded BGG correspondence. The functors $\LL$ and $\RR$ are given as follows. We will only need the formula for $\LL$ on $E$-modules with trivial differential; see \cite[\S 2]{BETate} for the formula in full generality. Given an $E$-module~$N$, the complex $\LL(N)$ has terms and differentials given by
$$
\LL(N)_j = \bigoplus_{a \in A}S(-a)  \otimes_k N_{(a; j)} \quad \text{ and } \quad s \otimes y \mapsto  (\sum_{i = 0}^n x_is \otimes e_iy).
$$
 Let $\om_E \ce \Hom_k(E, k) \cong E(-\sum_{i = 0}^n \deg(x_i); -n-1)$, i.e. the canonical module of $E$. Given an $S$-module $M$, the object $\RR(M) \in \DM(E)$ is $\bigoplus_{a \in A} M_a \otimes_k \om_E(-a; 0)$ equipped with the differential given by $m \otimes f \mapsto \sum_{i = 0}^n x_im\otimes e_if$. The package \verb|MultigradedBGG| implements these functors in \verb|Macaulay2|. Let us discuss some examples. 

\begin{example}
\label{ex:LL}
Applying the functor $\LL$ to a rank 1 free module should give the Koszul complex on the variables in $S$, up to a homological shift. We confirm this in the following example. We work over the Cox ring $S$ of a Hirzebruch surface of type 3, i.e. the ring $S = k[x_0, \dots, x_3]$ with $\Z^2$-grading given by $\deg(x_0) = \deg(x_2) = (1,0)$, $\deg(x_1) = (-3, 1)$, and $\deg(x_3) = (0,1)$. 
\begin{footnotesize}
\begin{verbatim}

i1 : needsPackage "MultigradedBGG"
i2 : needsPackage "NormalToricVarieties"
i3 : S = ring hirzebruchSurface 3;
i4 : E = dualRingToric S;

\end{verbatim}
\end{footnotesize}
The command \verb|dualRingToric S| gives the $\Z^3$-graded exterior algebra $E$ described above. 
\begin{footnotesize}
\begin{verbatim}
     
i5 : C = toricLL(E^1)

                 1                4                6                4                1
o5 = (QQ[x ..x ]) <-- (QQ[x ..x ]) <-- (QQ[x ..x ]) <-- (QQ[x ..x ]) <-- (QQ[x ..x ])
          0   3            0   3            0   3            0   3            0   3
                                                                               
      -4               -3               -2               -1               0

\end{verbatim}
\end{footnotesize}
\end{example}

\begin{example}
\label{ex:RR}
As in Example~\ref{ex:LL}, let $S = k[x_0, \dots, x_3]$ be the $\Z^2$-graded Cox ring of a Hirzebruch surface of type 3, and take $M = S / (x_0, x_1^2, x_2^2, x_3^2)$. Since $\dim_k M < \infty$, the differential $E$-module $\RR(M)$ has finite rank, and we can use the package \verb|MultigradedBGG| to compute $\RR(M)$ in its entirety. Since $\dim_k M = 8$, the rank of $\RR(M)$ is 8. 
\begin{footnotesize}
\begin{verbatim}

i1 : needsPackage "MultigradedBGG"
i2 : needsPackage "NormalToricVarieties"
i3 : S = ring hirzebruchSurface 3;
i4 : M = coker matrix{{x_0, x_1^2, x_2^2, x_3^2}};
i5 : (toricRR M).dd_0

o5 = {-1, 2, 4} | 0   0   0   0 0   0   0   0   |
     {-4, 3, 4} | e_1 0   0   0 0   0   0   0   |
     {-3, 3, 4} | 0   e_2 0   0 0   e_1 0   0   |
     {-3, 4, 4} | 0   0   e_3 0 e_2 0   e_1 0   |
     {-4, 4, 4} | 0   e_3 0   0 0   0   0   e_1 |
     {0, 2, 4}  | e_2 0   0   0 0   0   0   0   |
     {0, 3, 4}  | 0   0   0   0 0   e_3 0   e_2 |
     {-1, 3, 4} | e_3 0   0   0 0   0   0   0   |

\end{verbatim}
\end{footnotesize}
In the typical case where $\dim_k M = \infty$, the differential $E$-module $\RR(M)$ has infinite rank, and so one must pick a range of degrees in which to view $\RR(M)$. When $\dim_k M = \infty$, the method \verb|toricRR| gives, by default, the finite rank quotient of the differential $E$-module $\RR(M)$ given by the sum of 
 all $M_d \otimes_k \omega_E(-d; 0)$ such that $d = e + a\deg(x_i)$, where $e$ is a generating degree of $M$, $a \in \{0, 1\}$, and $0 \le i \le n$.
For instance, take $M = S / (x_0)$. Applying the method \verb|toricRR| to $M$ gives a rank 5 quotient of $\RR(M)$, since $M$ is generated in degree~0, $\dim_k M_{(0,0)} = \dim_k M_{(1,0)} = \dim_k M_{(-3,1)} = 1$, and $\dim_k M_{(0,1)} = 2$.
\begin{footnotesize}
\begin{verbatim}

i5 : M = coker matrix{{x_0}};
i6 : (toricRR M).dd_0
o6 = {-1, 2, 4} | 0   0 0 0 0 |
     {0, 2, 4}  | e_2 0 0 0 0 |
     {-4, 3, 4} | e_1 0 0 0 0 |
     {-1, 3, 4} | 0   0 0 0 0 |
     {-1, 3, 4} | e_3 0 0 0 0 |
     
\end{verbatim}
\end{footnotesize}
We may also apply the method \verb|toricRR| to a module $M$ along with a list $L$ of degrees in $A$: the output is the quotient of $\RR(M)$ given by $\bigoplus_{d \in L} M_d \otimes_k \om_E(-d ; 0)$. 
\begin{footnotesize}
\begin{verbatim}

i7 : L = {{0,0}, {1,0}, {-3, 1}, {0,1}, {2,0}};
i8 : (toricRR(M, L)).dd_0
o8 = {-1, 2, 4} | 0   0   0 0 0 0 |
     {0, 2, 4}  | e_2 0   0 0 0 0 |
     {-4, 3, 4} | e_1 0   0 0 0 0 |
     {-1, 3, 4} | 0   0   0 0 0 0 |
     {-1, 3, 4} | e_3 0   0 0 0 0 |
     {1, 2, 4}  | 0   e_2 0 0 0 0 |
\end{verbatim}
\end{footnotesize}
\end{example}

\begin{example}
Suppose $S = k[x_0, \dots, x_n]$ is $\Z$-graded, where $\deg(x_i) = d_i \ge 1$. In other words, $S$ is the Cox ring of the weighted projective space $X = \PP(d_0, \dots, d_n)$. Given a finitely generated $S$-module $M$, it follows from \cite[Theorem 3.7]{BETate} that the cohomology of the twists of the sheaf $\widetilde{M}$ on $X$ may be computed by taking a minimal free resolution of $\RR(M)$. In more detail: let $\e \co F \xra{\simeq} \RR(M)$ be the minimal free resolution of $\RR(M)$, and let $T(\widetilde{M})$ be the mapping cone of $\e$. The differential $E$-module $T(M)$ is called the \emph{Tate resolution} of the sheaf $\widetilde{M}$. By \cite[Theorem 3.7]{BETate}, we have: 
\begin{equation}
\label{eqn:coh}
H^j(X, \widetilde{M}(i)) \cong \Hom_E(k, T(\widetilde{M}))_{(i, -j)}.
\end{equation}
As a simple example, suppose $X = \PP(1, 1, 2)$ and $M = S/(x_0, x_1)$. We have 
\begin{equation}
\label{eqn:calc}
\dim_k H^j(X, \widetilde{M}(i)) = \begin{cases} 1, & j = 0 \text{ and }i \text{ even};\\ 0, & \text{else}.\end{cases}
\end{equation}
Let us check this via the isomorphism~\eqref{eqn:coh}. We have $\RR(M) = \bigoplus_{i \ge 0} \om_E(-2i, 0)$, with differential mapping $\om_E(-2i, 0)$ to $\om_E(-2i - 2, 0)$ via multiplication by $e_2$. The minimal free resolution of $\RR(M)$ is easily seen to be $F = \bigoplus_{i \ge 1} \om_E(2i, 1)$, with differential that sends $\om_E(2, 1)$ to zero and $\om_E(2i, 1)$ to $\om_E(2i - 2, 1)$ via multiplication by $e_2$ for $i \ge 2$. As an $E$-module, we have
$
T(\widetilde{M}) = F(0, -1) \oplus \RR(M) = \bigoplus_{i \in \Z} \om_E(2i, 0).
$
The calculation~\eqref{eqn:calc} thus follows from the isomorphism~\eqref{eqn:coh}.

One may use \verb|MultigradedBGG| to check all of this, but there are subtleties to be aware of. We first compute $\RR(M)$ (in a finite window):
\begin{footnotesize}
\begin{verbatim}

i1 : X = weightedProjectiveSpace {1,1,2};

i2 : S = ring X;

i3 : M = coker matrix{{x_0, x_1}};

i4 : D = toricRR(M, for i from 0 to 4 list i);

i5 : D.dd_0

o5 = {4, 3} | 0   0   0 |
     {6, 3} | e_2 0   0 |
     {8, 3} | 0   e_2 0 |
     
\end{verbatim}
\end{footnotesize}
Let us now resolve $\RR(M)$. While the function \verb|resDM| does not always yield a minimal free resolution, it does give a portion of the minimal free resolution of $\RR(M)$ in this case (the problem of non-minimality can be solved by slightly modifying Algorithm~\ref{alg:res2}, but we will not pursue this here). Since we have truncated $\RR(M)$, we have created additional homology that manifests itself as extraneous summands in the free resolution. 
\begin{footnotesize}
\begin{verbatim}

i6 : F = resDM(D, 3);

i7 : d = F.dd_0

o7 = {8, 3}  | 0 0 e_2 0   0   0   |
     {2, 2}  | 0 0 0   e_2 0   0   |
     {6, 3}  | 0 0 0   0   e_2 0   |
     {0, 2}  | 0 0 0   0   0   e_2 |
     {4, 3}  | 0 0 0   0   0   0   |
     {-2, 2} | 0 0 0   0   0   0   |
     
\end{verbatim}
\end{footnotesize}
The submatrix corresponding to the degree $(2,2), (0,2)$, and $(-2, 2)$ generators gives the summand we want; the rest corresponds to homology of our truncation of $\RR(M)$ that we wish to ignore.
\begin{footnotesize}
\begin{verbatim}

i9 : submatrix'(d,{0,2,4},{0,2,4})

o9 = {2, 2}  | 0 e_2 0   |
     {0, 2}  | 0 0   e_2 |
     {-2, 2} | 0 0   0   |
     
\end{verbatim}
\end{footnotesize}
Since these three summands have socle degree $(-2,-1), (-4,-1)$, and $(-6,-1)$, their presence implies that $H^0(X, \widetilde{M}(j)) = 1$ for $j = -2, -4, -6$. To compute the cohomology of additional negative twists, one needs to resolve a larger portion of $\RR(M)$. 

In general, to compute the cohomology of a sheaf $\widetilde{M}$ on weighted projective space using the methods in \verb|MultigradedBGG|, there are delicate matters to address. For instance, what is the most efficient window of $\RR(M)$ to use? How much of the (typically infinite) free resolution of $\RR(M)$ must one compute in order to be sure one has included all of the needed summands? How does one ensure that the extraneous summands resulting from the truncation of $\RR(M)$ do not affect the calculation? All of these matters will be addressed in future work of the second author with Daniel Erman. 
\end{example}

\section{Computing strongly linear strands of multigraded free resolutions}
\label{sec:strands}

Let $S$ be as in \S\ref{sec:BGG}, and assume $S$ is \emph{positively graded}, by which we mean that $S_0 = k$, and there exists a homomorphism $\theta\colon A \to \Z$ such that $\theta(\deg(x_i)) > 0$ for all $i$. For instance, if $S$ is the Cox ring of a smooth projective toric variety, then $S$ is positively graded by $\Cl(X)$ \cite[Example 5.2.3]{CLS}. Let $f \c F \to G$ be a homogeneous map of $A$-graded free $S$-modules. We say $f$ is \emph{strongly linear} if there exist bases of $F$ and $G$ with respect to which $f$ is a matrix whose entries are $k$-linear combinations of the variables.\footnote{There are several competing notions of linearity for maps over polynomial rings with nonstandard gradings; see \cite[\S 4]{BEsyzygies} for a detailed discussion.}  A complex $C$ of free $A$-graded $S$-modules is \emph{strongly linear} if
its differentials are such. Given a finitely generated graded $S$-module $M$ that is generated in a single degree, the \emph{strongly linear strand} of its minimal free resolution $F$ is, roughly speaking, the largest strongly linear subcomplex of $F$. More precisely, it is the unique maximal strongly linear subcomplex of $F$ that is a summand of $F$ (as an $S$-submodule, but not as a subcomplex) \cite[Definition 5.1]{BEstrands}. When $S$ is standard graded, this recovers the usual notion of the linear strand of a free resolution.

Let $M$ be an $A$-graded $S$-module that is generated in a single degree $d$. By \cite[Theorem 1.3(2)]{BEstrands}, there is a formula for the strongly linear strand $L$ of the minimal free resolution of $M$ in terms of the BGG correspondence. More precisely:
\begin{equation}
\label{linearformula}
L = \LL\left( \ker \left(  M_d \otimes_k \om_E(-d , 0) \xra{\del_{\RR(M)}} \RR(M) \right) \right).
\end{equation}
This generalizes a well-known result in the standard graded case: see \cite[Corollary 7.11]{geosyz}. As an application of our implementation of the multigraded BGG functors, we implement formula \eqref{linearformula} as the function \verb|stronglyLinearStrand| in \verb|MultigradedBGG|, allowing one to compute strongly linear strands (for modules generated in a single degree) in \verb|Macaulay2|. 

\begin{example}
As in Examples~\ref{ex:LL} and~\ref{ex:RR}, let $S = k[x_0, \dots, x_3]$ be the $\Z^2$-graded Cox ring of the Hirzebruch surface of type 3. The minimal free resolution of $M = S / (x_0, x_1^2)$, is, of course, the Koszul complex on $x_0$ and $x_1^2$. The strongly linear strand of this resolution is therefore the Koszul complex on $x_0$. Let us verify this with \verb|MultigradedBGG|:
\begin{footnotesize}
\begin{verbatim}

i1 : needsPackage "MultigradedBGG"
i2 : needsPackage "NormalToricVarieties"
i3 : S = ring hirzebruchSurface 3;
i4 : M = coker matrix{{x_0, x_1^2}};
i5 : (stronglyLinearStrand(M)).dd
                     1                          1
o5 = 0 : (QQ[x ..x ])  <----------- (QQ[x ..x ])  : 1
              0   3       | x_0 |        0   3

\end{verbatim}
\end{footnotesize}
As another example, let $S = k[x_0, \dots, x_4]$ be the Cox ring of the weighted projective space $X = \PP(1,1,1,2,2)$, so that $\deg(x_0) = \deg(x_1) = \deg(x_2) = 1$, and $\deg(x_3) = \deg(x_4) = 2$. We consider the curve $C \subseteq X$ defined by the $2 \times 2$ minors of the matrix $\begin{pmatrix} x_0 & x_1 & x_2^2 & x_3 \\ x_1 & x_2 & x_3 & x_4\end{pmatrix}$; let $I_C \subseteq S$ denote this ideal. The curve $C$ is isomorphic to $\PP^1$; indeed, it is the image of the closed immersion
$
\PP^1 \to X$ given by $[s : t] \mapsto [s^3 : s^2t : st^2 : st^5 : t^6]. 
$
A direct computation shows that $S/I_C$ is Cohen-Macaulay; let $M = \Ext^3_S(S/I_C, S(-7))$  be its canonical module. The strongly linear strand of the free resolution of $M$ is computed in \cite[Example 6.6]{BEstrands}; let us recover that computation using \verb|MultigradedBGG|.
\begin{footnotesize}
\begin{verbatim}

i6 : S = ring weightedProjectiveSpace {1,1,1,2,2};
i7 : I = minors(2, matrix{{x_0, x_1, x_2^2, x_3}, {x_1, x_2, x_3, x_4}});
i8 : M = Ext^3(module S/I, S^{{-7}});
i9 : (stronglyLinearStrand M).dd
                      3                                                    6
o9 =  0 : (QQ[x ..x ])  <------------------------------------- (QQ[x ..x ])  : 1
               0   4       {1} | x_0 0   x_3 0   -x_1 0    |        0   4
                           {1} | x_1 x_0 x_4 x_3 -x_2 -x_1 |
                           {1} | 0   x_1 0   x_4 0    -x_2 |

                      6                                       3
      1 : (QQ[x ..x ])  <------------------------ (QQ[x ..x ])  : 2
               0   4       {2} | x_3  0   x_1 |        0   4
                           {2} | x_4  0   x_2 |
                           {3} | -x_0 x_1 0   |
                           {3} | -x_1 x_2 0   |
                           {2} | 0    x_3 x_0 |
                           {2} | 0    x_4 x_1 |

\end{verbatim}
\end{footnotesize}
The differentials do not match exactly with those computed in \cite[Example 6.6]{BEstrands}, but the two complexes are isomorphic.
\end{example}

\bibliographystyle{amsalpha}
\bibliography{Bibliography}

\Addresses

\end{document}